\documentclass[12pt,leqno]{article}
\usepackage{amsthm,amsfonts,amssymb,amsmath,eufrak,oldgerm}
\usepackage{epsfig}
\numberwithin{equation}{section} 

\renewcommand\d{\partial}
\renewcommand\a{\alpha}

\renewcommand\b{\beta}

\renewcommand\o{\omega}
\newcommand\R{\mathbb R}
\newcommand\C{\mathbb C}

\def\eps{\varepsilon}
\def\e{\varepsilon}

\newcommand\br{\begin{remark}}
\newcommand\er{\end{remark}}
\newcommand\bp{\begin{pmatrix}}
\newcommand\ep{\end{pmatrix}}
\newcommand\be{\begin{equation}}
\newcommand\ee{\end{equation}}
\newcommand\ba{\begin{equation}\begin{aligned}}
\newcommand\ea{\end{aligned}\end{equation}}

\newcommand{\bap}{\begin{app}}
\newcommand{\eap}{\end{app}}
\newcommand{\begs}{\begin{exams}}
\newcommand{\eegs}{\end{exams}}
\newcommand{\beg}{\begin{example}}
\newcommand{\eeg}{\end{exaplem}}
\newcommand{\bpr}{\begin{proposition}}
\newcommand{\epr}{\end{proposition}}
\newcommand{\bt}{\begin{theorem}}
\newcommand{\et}{\end{theorem}}
\newcommand{\bc}{\begin{corollary}}
\newcommand{\ec}{\end{corollary}}
\newcommand{\bl}{\begin{lemma}}
\newcommand{\el}{\end{lemma}}
\newcommand{\bd}{\begin{definition}}
\newcommand{\ed}{\end{definition}}
\newcommand{\brs}{\begin{remarks}}
\newcommand{\ers}{\end{remarks}}

\newcommand{\CalG}{\mathcal{G}}

\newcommand{\RR}{{\mathbb R}}

\newcommand{\Id}{{\rm Id }}

\newcommand{\Res}{{\rm Residue}}

\newtheorem{theorem}{Theorem}[section]
\newtheorem{proposition}[theorem]{Proposition}
\newtheorem{corollary}[theorem]{Corollary}
\newtheorem{lemma}[theorem]{Lemma}
\newtheorem{definition}[theorem]{Definition}

\newtheorem{example}[theorem]{Example}
\newtheorem{remark}[theorem]{Remark}

\newtheorem{exams}[theorem]{Examples}

\newcommand\cB{{\cal  B}}

\newcommand\cG{{\cal  G}}
\newcommand\cK{{\cal  K}}

\newcommand\cN{{\cal  N}}

\newcommand\cF{{\cal  F}}

\newcommand\cM{{\mathcal M}}
\newcommand\cT{{\mathcal T}}

\newcommand\tG{{\tilde G}}

\title{
Conditional stability of unstable viscous shocks 
}

\author{\sc \small 
Kevin Zumbrun\thanks{Indiana University, Bloomington, IN 47405;
kzumbrun@indiana.edu:
Research of K.Z. was partially supported
under NSF grants no. DMS-0070765 and DMS-0300487.
 }}
\begin{document}

\maketitle

%%%%%%%%%%%%%%%%%%%%%%%%%%%%%%%%%%%%%%%%%%%%%%%%%%%%%%%%%%%%%%%%%%%%%%%%%%%%%%%%%%%%%%%%%%%%%

\begin{abstract}
Continuing a line of investigation initiated by Texier and Zumbrun 
on dynamics of viscous shock and
detonation waves, 
we show that a linearly unstable Lax-type viscous shock solution
of a semilinear strictly parabolic system of conservation laws
possesses a translation-invariant center stable manifold within which
it is nonlinearly orbitally stable with respect to small $L^1\cap H^2$ 
perturbatoins, converging time-asymptotically
to a translate of the unperturbed wave.
That is, for a shock with $p$ unstable eigenvalues, we establish
conditional stability on a codimension-$p$ manifold of initial data,
with sharp rates of decay in all $L^p$.
For $p=0$, we recover the result of unconditional stability obtained
by Howard, Mascia, and Zumbrun. 
\end{abstract}

%\clearpage
\tableofcontents
%\clearpage
%%%%%%%%%%%%%%%%%
\bigbreak
\section{Introduction}

In this paper, we continue a line of investigation opened
in \cite{TZ1,TZ2,TZ3,TZ4,SS,BeSZ}
going beyond simple stability analysis to study nontrivial dynamics, and
associated physical phenomena, of perturbed viscous shock waves in
the presence of {\it linear instability}.
The above-mentioned references concern Hopf bifurcation 
to time-periodic behavior associated with transition to
linear instability arising through the passage from stable
to unstable half-plane of a complex conjugate pair of eigenvalues
of the linearized operator about the wave.
%As discussed at a linear level in \cite{ZH,Z1},
%passage of a simple eigenvalue through zero may be associated either with
%wave-splitting or else appearance of multiple shock profiles, depending
%on the situation.
See also \cite{ZH,Z1} for discussion of (nonstandard, due to embedding
in essential spectrum) bifurcations
associated with the passage of a simple eigenvalue through zero. 

In the present work, in the spirit of \cite{GJLS,Li} and other
works outside the shock wave context, 
we consider the situation of a viscous shock substantially after the onset
of instability, i.e., {\it with one or more strictly unstable but 
no neutrally unstable eigenvalues}, and seek to describe the nearby
phase portrait in terms of invariant manifolds and behavior therein.
Specifically, for shock waves of systems of conservation laws with
artificial viscosity, we construct a center stable manifold and show that
the shock is conditionally (nonlinearly) stable with respect to this 
codimension $p$ set of initial data, where $p$ is the number of 
unstable eigenvalues.
As discussed for example in \cite{AMPZ,GZ}, such conditionally
stable shock waves can play an important role in asymptotic
behavior as metastable states.

Consider a viscous shock solution $u(x,t)=\bar u(x)$,
$\lim_{z\to \pm \infty}\bar u(z)=u_\pm$,
without loss of generality stationary,
of a semilinear parabolic system of conservation laws
\be\label{cons}
u_t+ f(u)_{x_j}= u_{xx},
\ee
$u,\, f\in \R^n$, $x\,,t\in \R$, under the basic assumptions:

\medbreak

(H0) $f\in C^{k+2}$, $k\ge 2$.
%NOTE: to get C^2 CS-manifold...

(H1) $A_\pm:=df(u_\pm)$ have simple, real, nonzero eigenvalues.
\medbreak

\noindent Linearizing \eqref{cons} about $\bar u$ yields linearized equations
\be\label{lineq}
u_t=Lu:=-(df(\bar u)u)_x -u_{xx},
\ee
for which the generator $L$ possesses \cite{He,Sat,ZH} 
both a translational zero-eigenvalue
and essential spectrum tangent at zero to the imaginary axis. 

The absence of a spectral gap between neutral (i.e., zero real part)
and stable (negative real part) spectra of $L$ prevents the usual
ODE-type decomposition of the flow near $\bar u$ into invariant
stable, center, and unstable manifolds.
The first result of this paper, by now little more
than 
%an observation, 
a remark,
is that we can still determine
center stable and unstable manifolds, and that these may
be chosen to respect the underlying translation-invariance of \eqref{cons}.
See \cite{TZ1} for closely related results on existence of
translational-invariant center manifolds.
As the needed ingredients do not seem to be found in one place,
we nonetheless 
for completeness
carry out the proof in full detail.

\bt\label{t:maincs}
Under assumptions (H0)--(H1), 
there exists in an $H^2$
neighborhood of the set of translates of $\bar u$
a codimension-$p$ translation invariant $C^k$ (with respect to $H^2$) 
center stable manifold $\cM_{cs}$,
tangent at $\bar u$ to the center stable subspace $\Sigma_{cs}$ of $L$, 
that is (locally) invariant under the forward time-evolution of \eqref{rcd}
and contains all solutions that remain bounded and sufficiently
close to a translate of $\bar u$ in forward time, where $p$ is the
(necessarily finite) number of unstable, i.e., positive real part,
eigenvalues of $L$.  
\et

Next, specializing a bit further, we add to (H0)--(H1) the additional
hypothesis that $\bar u$ be a {\it Lax-type shock}:
\medbreak

(H2) The dimensions of the unstable subspace of $df(u_-)$
and the stable subspace of $df(u_+)$ sum to $n+1$.
\medbreak

We assume further the following {\it spectral genericity} conditions.
\medbreak

(D1) $L$ has no nonzero imaginary eigenvalues. 

(D2) The orbit $\bar u(\cdot)$ is a transversal connection of
the associated standing wave equation $\bar u_x=f(\bar u)-f(u_-)$.

(D3) The associated inviscid shock $(u_-,u_+)$ is hyperbolically
stable, i.e.,
\be\label{liumajda}
\det(r_1^-,\dots, r_{p-1}^-, r_{p+1}^+, \dots , r_n^+, (u_+-u_-))\ne 0,
\ee
where $r_1^-, \dots r_{p-1}^-$ denote eigenvectors of $df(u_-)$
associated with negative eigenvalues and
$r_{p+1}^+, \dots r_{n}^+$ denote eigenvectors of $df(u_+)$
associated with positive eigenvalues.

\medbreak
\noindent As discussed in \cite{ZH,MaZ1}, (D2)--(D3) correspond in the
absence of a spectral gap to a generalized notion of simplicity of the 
embedded eigenvalue $\lambda=0$ of $L$.
Thus, (D1)--(D3) together correspond to the assumption that there are
no additional (usual or generalized) eigenvalues on the imaginary
axis other than the transational eigenvalue at $\lambda=0$;
that is, the shock is not in transition between different degrees of
stability, but has stability properties that are insensitive to 
small variations in parameters.

With these assumptions, we obtain our second and main result
characterizing stability properties of $\bar u$.
In the case $p=0$, this reduces to the nonlinear orbital stability
result established in \cite{ZH,MaZ1,MaZ2,MaZ3,Z2,Z3,Z4}.

\bt\label{t:mainstab}
Under (H0)--(H2) and (D1)--(D3), $\bar u$ is nonlinearly
orbitally stable under sufficiently small perturbations
in $L^1\cap H^4$ lying on the codimension $p$ center stable 
manifold $\cM_{cs}$ of $\bar u$ and its translates, 
where $p$ is the number of unstable eigenvalues of $L$,
in the sense that, for some $\alpha(\cdot)$,
\ba\label{bounds}
|u(x, t)-\bar u(x-\alpha(t))|_{L^p}&\le
C(1+t)^{-\frac{1}{2}(1-\frac{1}{p})}
|u(x,0)-\bar u(x)|_{L^1\cap H^4},
\\
|u(x, t)-\bar u(x-\alpha(t))|_{H^4}&\le
C(1+t)^{-\frac{1}{4}} |u(x,0)-\bar u(x)|_{L^1\cap H^4},\\
\dot \alpha(t) &\le C(1+t)^{-\frac{1}{2}} |u(x,0)-\bar u(x)|_{L^1\cap H^4},
\\
\alpha(t) &\le C |u(x,0)-\bar u(x)|_{L^1\cap H^4}.
\ea
Moreover, it is orbitally unstable with respect to small $H^2$ perturbations
not lying in $\cM$, in the sense that the corresponding solution leaves
a fixed-radius neighborhood of the set of translates of $\bar u$ in
finite time.
\et

%TODO: reinstate?
%NO, not needed either I think..
\subsection{Discussion and open problems}\label{discuss}

It is easily checked that the results of this paper go through 
for general semilinear parabolic term $Bu_{xx}$, $B$ constant, 
under the standard assumptions of \cite{ZH}.
An interesting open problem is to extend these results to the
quasilinear and or partially parabolic (``real viscosity'') case.
We note that the only issue here is to establish existence of
the center stable manifold, as the proof of stability goes through
essentially unchanged, incorporating the necessary modifications
detailed in \cite{Z2,Z3} to deal with the quasilinear/partially 
parabolic case.
See also related discussion in \cite{Li}.

Another interesting problem would be to extend our conditional
stability result to the case of nonclassical under- or overcompressive
shocks using pointwise estimates as in \cite{HZ,RZ}; see Remark \ref{ucpt}.
%TODO: more on this aspect?

Finally, we mention the problem of determining conditional stability
of a planar standing shock $u(x,t)\equiv \bar u(x_1)$ of a
multidimensional system of conservation laws
$$
u_t + \sum_j f_j(u)_{x_j}= \Delta_x u,
$$
which likewise (by the multidimensional arguments of \cite{Z1,Z2,Z3})
reduces to construction of a center stable manifold, in this case
involving an infinite-dimensional unstable subspace corresponding
to essential spectra of the linearized operator $L$ about the wave.

\medskip

{\bf Plan of the paper.}
In Section \ref{s:ode} we give for completeness a
particularly concise proof of the center stable manifold
theorem for ODE.
In Section \ref{s:pde}, loosely following \cite{TZ1},
we show how to extend this to semilinear
parabolic PDE, while preserving the key property of translation-invariance.
Finally, in Section \ref{s:cond}, we establish conditional stability
by a modification of the arguments of \cite{Z4,MaZ2,MaZ3} in the stable 
($p=0$) case.

%%%%%%%%%%%%%%%%%%%

\section{Center Stable Manifold for ODE}\label{s:ode}

The center stable manifold construction in the PDE case follows
closely the construction for finite-dimensional ODE, 
which we therefore recall here for completeness;
see also \cite{B,VI}.
%following a modified version of \cite{B,VI}.
%We follow a streamlined version of the center manifold treatments 
%in \cite{B,VI} which, for its brevity, may be of some interest independently.
%We spend particular care on the issue of smoothness, 
%since it is critical for the eventual
%application to stability that the constructed center stable manifold be 
%at least $C^2$.
%
Consider an ODE 
\be\label{nonline}
u'=f(u),\qquad f\in C^1,
\ee
and an equilibrium $f(u_*)=0$, with associated linearized equation 
\be\label{line}
v'=Av,
\qquad A:=df(u_*).
\ee
Associated with $A$, define the center stable subspace $\Sigma_{cs}$
as the direct sum of all eigenspaces of $A$ associated to
neutral or stable eigenvalues, i.e., eigenvalues with zero or positive real
part.  Likewise, define the unstable
subspace $\Sigma_u$ as the direct
sum of eigenspaces associated to unstable eigenvalues,
i.e., eigenvalues with strictly positive
real part, so that $\RR^n= \Sigma_{cs} \oplus \Sigma_u$.

Defining the associated (total) eigenprojections  $\Pi_{cs}$
and $\Pi_u$ as the sum of all eigenprojections 
associated with neutral--stable and unstable eigenvalues, respectively,
we have, 
either by reduction to Jordan form
or direct estimation using the inverse Laplace transform/resolvent estimates,
bounds
\ba\label{expbds}
|e^{At}\Pi_{cs}|&\le C(\eta)e^{\theta t} \qquad t\ge 0,\\
|e^{At}\Pi_u|&\le C(\theta) e^{-\eta |t|}, \quad t\le 0,\\
\ea
for any $\eta>0$ strictly smaller than the minimum 
of the real parts of unstable eigenvalues and
$\theta>0$ arbitrarily small.

\bpr [Center Stable Manifold Theorem for ODE] \label{t:centerman}
For 
%NOTE: typesetting choice...
%$f\in C^{k+1}$, $1\le k < \infty$, there exists local to $u_*$
$f$ in $C^{k+1}$, $1\le k < \infty$, there exists local to $u_*$
a $C^k$ center stable manifold $\cM_{cs}$,
tangent at $u_*$ to $\Sigma_{cs}$, expressible
in coordinates $w:=u-u_*$ as a $C^k$ graph 
$\Phi_{cs}:\Sigma_{cs}\to \Sigma_{cs}\oplus \Sigma_u$,
that is (locally) invariant 
%NOTE: in forward time, hence also backward, by dimensionality argument...
%this last part doesn't work for infinite-dim.!
under the flow of \eqref{nonline}
and contains all solutions that remain bounded and sufficiently
close to $u_*$ in forward time.
In general it is not unique.
\epr

\subsection{Frechet differentiability of substitution operators}\label{sub}
Standard invariant (e.g., stable, center, center stable) manifold 
constructions proceed by fixed point iteration 
on various weighted $L^\infty$ spaces 
\be\label{etadef}
\|f\|_\eta:=\sup_{t\ge t_0}e^{\eta (t-t_0)}|f(t)|,
\quad \cB_\eta:=\{f: \, \|f\|_\eta<\infty\},
\ee
with $\eta$ positive for stable manifolds and negative for center
or center stable manifolds.
%On the latter spaces, the key issue for us 
As described in \cite{B,VI}, Frechet differentiability of the
associated fixed-point mapping
(hence eventual smoothness of the resulting manifold)
hinges on Frechet differentiability with respect to spaces $\cB_\eta$
of the special class of {\it substitution operators}, 
%TODO: mention later...
%NOTE: reason is that our integral mappings are all compositions of
%bounded linear operators on $\cB_\eta$ and evaluation maps.
%Thus, diff. of $T$ amounts (by composition rule) just to diff. of component
%evaluation map...!
defined for $g:\RR^n\to \RR^n$ as 
\be\label{evaldef}
G(f)(t):=g(f(t)).
\ee

Let $C^k_b$ denote the Banach space of $C^k$ functions $g:\RR^n\to \RR^n$
with $|d^jg|$ uniformly bounded for $0\le j\le k$,
with associated norm
$$
\|g\|_{C^k_b}:=\sum_{0\le j\le k} \sup_{\RR^n}|d^jg|.
$$

\bl\label{poseta}
For $\eta\ge 0$, $k\ge 1$, if $g\in C^k_b(\RR^n\to \RR^n)$, then
$G$ is $C^k$ from $\cB_\eta\to \cB_\eta$, with $d^kGf(t)=(d^kg)(f(t))$.
\el

\begin{proof}
See Appendix \ref{extraproofs}.
\end{proof}

\bl\label{negfrechet}
Let $g\in C^{k+1}_b$ and $0<-\eta'< -\eta/k$.  
Then, the substitution operator $G$ is $C^k$ from $\cB_{\eta'}\to
\cB_\eta$, with $d^kGf(t)=(d^kg)(f(t))$.
\el

\begin{proof}
More generally, the result of Lemma \ref{negfrechet} holds
for any $g\in C^{k+\alpha}_b$, $0<\alpha \le 1$, in the
sense that $|d^kg(x+h)-d^kg(x)|\le C|h|^\alpha$ for some 
uniform $C>0$, as may be seen by rewriting the $(k-1)$th-order 
Taylor remainder formula
$
g(x+h)-\cT_{k-1}g(x,h)=
\Big(\int_{0}^{1}d^{k}g(x+ \theta h)
\frac{(1-\theta)^{k-1}}{(k-1)!}d\theta\Big) h^{k}
$
as
$$
g(x+h)-\cT_{k}g(x,h)=
\Big(\int_{0}^{1}\big(d^{k}g(x+ \theta h)-d^kg(x)\big)
\frac{(1-\theta)^{k-1}}{(k-1)!}d\theta\Big) h^{k},
$$
where $\cT_k g(x,h)$ denotes the $k$th-order Taylor expansion of
$g$ about $x$ evaluated at $x+h$,
then using the assumed uniform bound on $|d^kg(x+h)-d^kg(x)|$ to
obtain $|g(x+h)-\cT_{k}g(x,h)|\le C|h|^{k+\alpha}$, and thus
$$
|g(x+h)-\cT_{k}g(x,h)|_{-(k+\alpha) \eta'}\le C|h|_{\eta'}^{k+\alpha}.
$$
Without loss of generality taking $\alpha$ sufficiently small,
this yields the result.
A similar estimate yields continuity of the $k$th Frechet derivative.
\end{proof}

\subsection{Smooth dependence 
of fixed-point solutions}\label{dependence}

We next present two general results on smooth dependence of fixed
point solutions.
%point solutions, the first completely standard and the second 
%concerning mappings that are smooth between different norms.
Let $T(x,y)$ be continuous in $x$, $y$ and contractive in $y$, 
$T:\cB_1\times \cB_2 \to \cB_2$ for Banach spaces $\cB_1$ and $\cB_2$,
defining a fixed point map $y(x)$, $y:\cB_1\to \cB_2$,
continuous in the parameter $x$, such that 
$y(x)=T(x,y(x))$.
Then, we have the following standard result. 

\bl\label{fixeddep}
If $T$ is Lipschitz in $(x,y)$, 
then $y$ is Lipschitz from $\cB_1\to \cB_2$.
If $T$ is $C^k$ (Frechet sense), $k\ge 1$, in $(x,y)$, then $y$
is $C^k$ from $\cB_1\to \cB_2$, with
\be\label{derform}
(dy/dx)(x_0)= (\Id-T_y)^{-1}T_x (x_0, y(x_0))
\ee
and higher derivatives $(d^j y/dx^j)(x_0)$, $1\le j\le k$ obtained
by formal differentiation of \eqref{derform}, substituting for 
lower derivatives wherever they appear.
\el

\begin{proof}
See Appendix \ref{extraproofs}.
\end{proof}

The next Lemma shows how we can recover $C^1$ dependence of
fixed point solutions in the case, as in Lemma \ref{negfrechet}, 
that $T$ is differentiable only from a stronger to a weaker space.
We discuss higher derivatives later where they appear,
since they involve specific chains of successively weaker spaces
that are not convenient for statement as a general theorem.
Let $\cB_2'\subset \cB_2$, with $\|\cdot\|_{\cB_2'}\ge \|\cdot\|_{\cB_2}$,
and $T(x,y)$ be a map $\cB_1\times \cB_2'\to \cB_2'$
that is Lipschitz continuous in $(x,y)$ and contractive 
(with respect to $\|\cdot\|_{\cB_2'}$) in $y$.
Denote by $y(x):\cB_1\to \cB_2'$ the unique 
Lipschitz fixed-point solution defined by $y(x)=T(x,y(x))$.

\bl\label{diffdep}
If (i) $T$ is continuously differentiable from $\cB_1\times \cB_2'\to \cB_2$,
and (ii) $T_y$ extends to a bounded linear operator from
$\cB_2\to \cB_2$, continuous in operator norm with respect to $(x,y)$, 
with $|T_y|_{\cB_2}<1$, then $y$ is continously
differentiable from $\cB_1\to \cB_2$, 
with $(dy/dx)(x_0)= (\Id-T_y)^{-1}T_x (x_0, y(x_0))$.
\el

\begin{proof}
%TODO.
By Taylor's Theorem and
$\|y(x_2)-y(x_1)\|_{\cB_2'}\le L\|x_2-x_1\|_{\cB_1}$, 
we have
\ba\label{cal}
y(x_2)-y(x_1)&=T_x(x_2-x_1) + T_y(y_2-y_1)\\
&\quad  +o(\|x_2-x_1\|_{\cB_1}
+\|y(x_2)-y(x_1)\|_{\cB_2'})\\
&=
T_x(x_2-x_1) + T_y(y_2-y_1) +o(\|x_2-x_1\|_{\cB_1}),
\ea
where the $o(\|x_2-x_1\|_{\cB_1})$ term is measured in the
weaker $\|\cdot\|_{\cB_2}$ norm.
Observing by $|T_y|_{\cB_2\to \cB_2}<1$ and Neumann series inversion 
that $(\Id -T_y)$ considered as an operator from $\cB_2\to \cB_2$
is invertible with uniformly bounded inverse
$$
|(\Id -T_y)^{-1}|_{\cB_2\to \cB_2}\le (1-|T_y|_{\cB_2\to \cB_2})^{-1},
$$
we may solve \eqref{cal} to obtain
$$
y(x_2)-y(x_1)= (\Id-T_y)^{-1}T_x (x_2-x_1)+ o(\|x_2-x_1\|_{\cB_1}),
$$
yielding differentiability as claimed, with
$(dy/dx)= (\Id-T_y)^{-1}T_x $ continuous by the assumed continuity
of $T_x$ and $T_y$ as operators from $\cB_1\to \cB_1$
and $\cB_2\to \cB_2$.
\end{proof}

\subsection{Global Center Stable Manifold construction}\label{cm}

%To prove Theorem \ref{t:centerman}, we first establish a global version for
We now establish a global version of \ref{t:centerman} for
small Lipschitz nonlinearity. 
As in Section \ref{dependence}, denote by $C^k_b$ the space 
of $C^k$ functions that are uniformly bounded in up to $k$ derivatives
and consider for a fixed, constant matrix $A$ and an arbitrary nonlinearity
$N$ such that $N(t,0)\equiv 0$, $N_w(t,0)\equiv 0$ the ODE
\be\label{eqcm}
w'=Aw+N(t,w).
\ee

\bpr\label{pr:globalcm}
For $N\in C^{k+1}_b$ with Lipschitz
constant $\eps>0$ sufficiently small, \eqref{eqcm}
has a unique $C^k$ invariant manifold $\cM_{cs}$ 
tangent at $w=0$ to the center stable subspace $\Sigma_{cs}$ of $A$,
consisting of the union of all orbits whose solutions grow
at sufficiently slow exponential rate $|w(t)|\le C e^{\tilde \theta |t|}$
in positive time,
for any fixed $\theta<\tilde \theta<\eta$.
%NOTE: $C>0$ here depending on the solution.
\epr

\begin{proof}
Applying projections $\Pi_j$, $j=cs,u$ to \eqref{eqcm},
we obtain using the variation of constants formula equations
$$
\Pi_j w(t)=
e^{A(t-t_{0,j})}\Pi_j w(t_{0,j}) + \int_{t_{0,j}}^t e^{A(t-s)}\Pi_j  N(s, w(s))\,ds,
$$
$j=cs,u$, so long as the solution $w$ exists,
with $t_{0,j}$ arbitrary.
Assuming growth of at most $|w(t)|\le Ce^{\tilde \theta t}$ in positive 
time, we find using \eqref{expbds} and the bound $|N(w)|\le \eps |w|$
coming from $N(0)=0$ and the assumed Lipschitz bound on $N$,
that as $t_{0,u}\to +\infty$, the first term
$e^{A(t-t_{0,s})}\Pi_u w(t_{0,j})$ converges to zero while
the second, integral term
converges to $\int_{t}^{+\infty} e^{A(t-s)}\Pi_u  N(s, w(s))\,ds$,
so that, choosing $t_{0,cs}=0$, we have
$$
\begin{aligned}
\Pi_{cs} w(t)&=
e^{At}\Pi_{cs} w(0) + \int_{0}^t e^{A(t-s)}\Pi_{cs}  N(s, w(s))\,ds,\\
\Pi_u w(t)&=  -\int_t^{+\infty} e^{A(t-s)}\Pi_u  N(s, w(s))\,ds.
\end{aligned}
$$

Summing, we obtain for $w_{cs}:=\Pi_{cs} w(0)$ the fixed-point representation
\ba\label{Tcrep}
w(t)=T(w_{cs}, w)&:=
e^{At}w_{cs} 
+ \int_{0}^t e^{A(t-s)}\Pi_{cs}  N(s, w(s))\,ds\\
&\quad
-\int_{t}^{+\infty}e^{A(t-s)}\Pi_{u}  N(s, w(s))\,ds,
\ea
valid for solutions growing 
%exponentially 
at rate at most $|w(t)|\le C e^{\tilde \theta t}$ 
%in positive time $t\ge 0$. 
for $t\ge 0$. 

Define now the negatively-weighted sup norm 
$$
\|f\|_{-\tilde \theta}:=\sup_{t\ge 0}e^{-\tilde \theta t}|f(t)|,
$$
noting that $ |f(t)|\le e^{\tilde \theta t} \|f\|_{-\tilde \theta}$ 
for all $t\ge 0$.
We first show that, for $|w_{cs}|$ and $\delta, \eps>0$ sufficiently small, 
the integral operator $T$ is Lipschitz in $w_{cs}$ and
contractive in $w$ from the $\|\cdot\|_{-\tilde \theta}$-ball 
$B(0,\delta)$ to itself.

Using \eqref{expbds}, $|N(t, w(t))|\le \epsilon |w(t)|$, and
$ |w(t)|\le e^{\tilde \theta |t|} \|w\|_{-\tilde \theta}$,
we obtain
\ba\nonumber
|T(w)(t)|&\le C e^{\theta |t|}|w_{cs}|
+C\epsilon\|w\|_{-\tilde \theta }
\Big( \int_{0}^t e^{  \theta|t-s|} e^{\tilde \theta |s|}\, ds
+\int_t^{+\infty} e^{\eta(t-s)} e^{\tilde \theta |s|}\, ds\Big),\\
\ea
hence, using $\eta\pm \tilde \theta >0$
and taking $C_1\epsilon <1/2$ and
$C|w_0|\le \delta/2$, that
\ba \nonumber
\|T(w)(t)\|_{-\tilde \theta}&\le C|w_{cs}|
+
C\epsilon\|w\|_{-\tilde \theta }
\Big(
\int_{0}^t e^{  \theta|t-s|} e^{\tilde \theta (|s|-|t|)}\, ds\\
&\quad
+\int_t^{+\infty} e^{\eta(t-s)} e^{\tilde \theta (|s|-|t|)}\, ds\Big)
%\\ &
\le C |w_{cs}| + C_1 \epsilon\|w\|_{-\tilde \theta } < \delta.
\ea
%TODO: do this part by contraction plus est. on first term?
%OR, put as general exercise: given good bound on $T(0)$ plus
%contractivity, show, etc. etc. (YES, this I think... enough to
%check contraction plus inhomogeneous term...)
Similarly, we find that
\ba\nonumber
\|T(w_1)-T(w_2)\|_{-\tilde \theta}
&\le C\epsilon\|w_1-w_2\|_{-\tilde \theta }
\Big(
\int_{0}^t e^{  \theta|t-s|} e^{\tilde \theta |s|}\, ds 
%\\ &\quad 
+\int_t^{+\infty} e^{\eta(t-s)} e^{\tilde \theta |s|}\, ds\Big)
\\ &
\le C_1\eps \|w_1-w_2\|_{-\tilde \theta } < (1/2)\|w_1-w_2\|_{-\tilde \theta },
\\
\ea
yielding contraction on $B(0,\delta)$ 
and thus existence of a unique fixed point $w=w(w_{cs})$.
A similar estimate shows that $T$ is Lipschitz in $w_{cs}$,
so that $w(\cdot)$ is Lipshitz from $\Sigma_{cs}$ to $\cB_{-\tilde \theta}$.

We next investigate smoothness of $w(\cdot)$.
Note that $T$ decomposes into the sum of a 
bounded, hence $C^\infty$, linear map $w_{cs} \to e^{At}w_{cs}$ 
from $\Sigma_{cs}\to \cB_{-\tilde \theta}$ and the composition
$\cK\cdot \cN$ of a bounded (hence $C^\infty$) linear map
$$
\cK(f):=
\int_{0}^t e^{A(t-s)}\Pi_{cs}  f(s)\,ds
-\int_{t}^{+\infty}e^{A(t-s)}\Pi_{u}  f(s)\,ds
$$
from $\cB_{-\tilde \theta} \to \cB_{-\tilde \theta}$
and a substitution operator $\cN(w)(s):=N(w(s))$
with $N\in C^{k+1}_b$
that by  Lemma \ref{negfrechet} is $C^k$ from 
$\cB_{-\tilde \theta/(k+1)}\to \cB_{-\tilde \theta}$.
Moreover, the first derivative $d\cN(w)(s)=dN(w(s))$,
by $|dN|\le \eps$, extends as a bounded linear operator 
from $\cB_{-\theta'}\to \cB_{-\theta'}$, any $\theta'>0$,
with $|dN|_{\cB_{-\theta'}}\le \eps$, whence $T_{w}$
extends as a bounded linear operator
from $\cB_{-\theta'}\to \cB_{-\theta'}$, that for any 
given $\theta'>0$, in particular $\tilde \theta$, is {\it contractive},
$|T_w|_{\cB_{-\theta'}}<1$, for $\eps>0$ sufficiently small, 
independent of $w$ and $w_{cs}$.
Applying Lemma \ref{diffdep}, we find that $w(w_{cs})$ is $C^1$ from
$\Sigma_{cs}\to \cB_{-\theta'}$ for any $\theta'>0$ and $\eps>0$ 
sufficiently small, with
\be\label{weq}
dw=(\Id-T_w)^{-1}T_{w_{cs}}.
\ee

Differentiating \eqref{weq} using the chain, 
product, and inverse linear operator derivative formulae, 
validated by the fact that $(\Id-T_w)^{-1}$ is
uniformly bounded (since $|T_w|\le \gamma<1$) from $\cB_{-\theta'}\to
\cB_{-\theta'}$ for all $0\le \theta_0<\theta'<\eta_0<\eta$,
and the observation that higher (mixed) partial derivatives of $T$
exist and are continuous from $\cB_{-\theta'/(k+1)}\to \cB_{-\theta'}$
on the same range, we find that $w(w_{cs})$ is $C^k$ from
$\cB_{-\tilde \theta/C}\to \cB_{-\tilde \theta}$ for $C>0$ sufficiently
large, so long as $\theta_0>0$ is chosen $\le \tilde \theta/C$ and
$\eps>0$ is taken sufficiently small.

Finally, defining
\ba\label{Phidefc}
\Phi(w_{cs})&:= \Pi_{u}w(w_s)|_{t=0} =
 -\int_{0}^{+\infty}e^{A(t-s)}\Pi_{u}  N(s, w(s))\,ds,
\ea
we obtain as claimed a $C^k$ function from $\Sigma_{cs}\to \Sigma_{u}$,
whose graph is the invariant manifold of orbits growing at exponential
rate $|w(t)|\le Ce^{\tilde \theta t}$ in forward time.
From the latter characterization, we obtain evidently
invariance in forward and backward time.
By uniqueness of fixed point solutions, we have $w(0)=0$ and
thus $\Phi_{cs}(0)=0$.
Finally, differentiating \eqref{Phidefc} with respect to $w_s$, we obtain
by Lemma \ref{negfrechet}
$$
d\Phi(0) =
 -\int_{0}^{+\infty}e^{A(t-s)}\Pi_{u}  N_w(s, 0)(dw/dw_s)(s)\,ds=0 ,
$$
since $N_w(0)\equiv 0$, yielding tangency as claimed.
\end{proof}

\subsection{Local construction: proof of Theorem \ref{t:centerman}}\label{lcm}

We can reduce the general situation, locally, to the case
described in the global result by the following truncation procedure.
Consider a general nonlinearity $N(t,w)$.
Introducing a $C^\infty$ cutoff function
$$
\rho(x) =\begin{cases}
1  & | x | \leq 1, \\ 
0 &  | x | \geq 2,\end{cases}
$$
define $N^\eps(t,w) := \rho(|w|/\eps)N(t,w)$.

\bl\label{truncbds}
Let $N\in C^{k+1}$, $k\ge 1$, and 
$N(t,0)\equiv 0$, $\partial_w N(t,0)\equiv 0$. Then, $N^\eps\in C^{k+1}_b$,
$N^\eps\equiv N$ for $|u|\le \eps$, and the Lipschitz constant 
for $N^\eps$ with respect to $w$ is uniformly bounded by $C\eps$ 
for all $0<\eps\le \eps_0$, where 
\be\label{Nebd}
 C=2\Big(1+ \max |r|\rho'(|r|)\Big) 
\max_{|w|\le 2\eps_0}|\partial_w^2N(t,w)|.
\ee
\el

\begin{proof}
See Appendix \ref{extraproofs}
\end{proof}

%\begin{proof}
%The Lipshitz constant is bounded by $\max |\partial_w N^\eps|$, where
%$$
%\begin{aligned}
%|\partial_w N^\eps|&= |(\rho^\eps)' N+ \rho^\eps \partial_w N|\\
%&= |(|w|/\eps)\rho'(|w|/\eps) (N(t,w)/|w|)+ \rho(|w|/\eps) \partial_w N(t,w)|\\
%&\le 2\eps \Big(
%\max |r|\rho'(|r|) \max_{|w|\le 2\eps_0}\frac{|N(t,w)|}{|w|^2}
%+ \max_{|w|\le 2\eps_0}\frac{|\partial_w N(t,w)|}{|w|}\Big)\\
%&\le 2\eps \Big(
%\max |r|\rho'(|r|) \max_{|w|\le 2\eps_0}
%+ 1 \Big) \max_{|w|\le 2\eps_0}|\partial_w ^2N(t,w)|,
%\end{aligned}
%$$
%the final inequality following by $N(t,0)\equiv 0$, 
%$\partial_w N(t,0)\equiv 0$ and
%the Integral Mean Value Theorem, or first-order Taylor
%remainder formula.
%\end{proof}

\begin{proof}[Proof of Proposition \ref{t:centerman}]
Defining $w:=u-u_*$, we obtain the nonlinear perturbation equation
\be\label{orig}
w'=Aw + N(w),
\ee
with $A$ constant and $N\in C^{k+1}$ satisfying $N(0)=0$, $dN(0)=0$.
Applying the truncation procedure, we obtain a modified equation
\be\label{mod}
w'=Aw + N^\eps(w)
\ee
for which $N^\eps\in C^{k+1}_b$ with arbitrarily small Lipschitz norm
$\eps>0$ and $N\equiv N^\eps$ within a neighborhood $B(0,\eps)$ 
of $w=0$, i.e., with identical local flow.
Applying Proposition \ref{pr:globalcm}, we obtain a global center
stable manifold for
\eqref{mod}, which is therefore a local center
stable manifold for \eqref{orig}.
Noting that solutions that stay uniformly bounded and close to the
equilibrium for positive time
are also bounded solutions of the truncated equations,
we find that all such belong to the constructed center stable manifold.
\end{proof}

\br
\textup{
The inclusion of $t$-dependence of $N$ in \eqref{eqcm} is not
needed for the present application \eqref{orig}, 
but allows also the treatment of time-periodic solutions 
after Floquet transformation to constant-coefficient linear part.
}
\er

\section{Center Stable Manifold for PDE}\label{s:pde}

We now turn to the case of a general semilinear parabolic equation 
\begin{equation}\label{rcd}
u_t = \cF(u):= h(u, u_x)+ u_{xx},
\qquad u, h\in \R^n,
\end{equation}
with steady state $u(x,t)\equiv \bar u(x)$ and  associated group invariance
\begin{equation}\label{transtilde}
\Psi_\alpha :\quad
\Psi_\alpha (u)(x,t):=
u(x + \alpha, t),
\end{equation}
and construct a local translation-invariant center stable manifold
in the vicinity of $\bar u$.

Our approach follows closely that used to construct
translation-invariant center manifolds in \cite{TZ1},
by introduction of a reduced flow on the quotient space
induced by group equivalence.
However, we coordinatize differently, using orthogonal projection
rather than eigenprojections, to avoid the difficulty
(as in the ultimate application to shock waves) that 
the zero eigenvale associated with translation-invariance
may be embedded in essential spectrum of the linearized
operator about the wave, as a consequence of which there
may not exist a well-defined zero eigenprojection with
respect to $H^s$ (see \cite{ZH} for further discussion).

We make the following assumptions, in practice typically satisfied \cite{TZ1}.

\medbreak
(A0) $h \in C^{k+1}$, $k\ge 2$.
%TODO: sort this out with earlier notation...

(A1) The linearized operator $L=\frac{\d \cF}{\d u}(\bar u)$ about
$\bar u$ has $p$ unstable (positive real part)
eigenvalues, with the rest of its spectrum  of nonpositive real part.
%Good (at least neutral) essential spectrum
%(Note, weaker than in \cite{TZ1}).
%TODO: explain. (just means assumption on $df(u_\pm)$,
%as discussed in \cite{TZ1}.)

(A2) $|\partial_x^j \bar u (x)|\le Ce^{-\theta |x|}$, $\theta>0$,
for $1\le j\le k+2$.

\medbreak

%TODO: be sure to mention it is a semigroup now and not a group
%action, so get invariance only in positive (forward) time direction...
%In particular, doesn't necessarily exist for negative times
%(not considered in the proof).

\bpr [Center Stable Manifold Thm. for PDE] \label{t:cspde}
Under assumptions (A0)--(A2), 
there exists 
%local to the set of translates of $\bar u$
in an $H^2$ neighborhood of the set of translates of $\bar u$
a translation invariant $C^k$ (with respect to $H^2$) 
%a translation invariant $C^k$ 
center stable manifold $\cM_{cs}$,
tangent at $\bar u$ to the center stable subspace $\Sigma_{cs}$ of $L$, 
that is (locally) invariant under the forward time-evolution of \eqref{rcd}
and contains all solutions that remain bounded and sufficiently
close to a translate of $\bar u$ in forward time.  In general it is not unique.
\epr

\subsection{Reduced equations}\label{redeq}

Differentiating with respect to $\alpha$ the relation
$
0\equiv \partial_t(\Psi^\alpha(\bar u)) = \cF(\Psi^\alpha(\bar u)) $,
we recover the standard fact that
$$
\phi := \frac{d \Psi_\a(\bar u)} {d \a}_{| \a = 0}= 
\frac{\partial \bar u}{\partial x} 
$$
is an $L^2$ zero eigenfunction of $L$, by the assumed decay of $\bar u_x$.

Define orthogonal projections
\begin{equation}
\label{proj}
\Pi_2 :=
\frac{ \phi \, \langle \phi, \cdot\rangle}{|\phi|_{L^2}^2}, 
\qquad \Pi_1 := \Id - \Pi_2,
\end{equation}
onto the range of right zero-eigenfunction $\phi:=(\partial/\partial x)\bar u$
of $L$ and its orthogonal complement $\phi^\perp$ in $L^2$, where 
$\langle \cdot, \cdot \rangle$ denotes standard $L^2$ inner product.

\begin{lemma}\label{Pi}
Under the assumed regularity $h\in C^{k+1}$, $k\ge 2$,
$\Pi_j$, $j=1,2$ are bounded as operators from $H^s$ to itself
for $0\le s \le k+2$.
\end{lemma}

\begin{proof}
Immediate, by the assumed decay of $\phi=\bar u_x$ and derivatives.
\end{proof}

Introducing the shifted perturbation variable
\be\label{shift}
v(x,t):= u(x+\alpha(t),t)-\bar u(x)
\ee
similarly as in \cite{Z4,MaZ2,TZ1}, 
we obtain the nonlinear perturbation equation
\begin{equation}
\label{epspert}
\d_t v= Lv + \CalG(v) -\d_t \alpha (\phi+ \partial_x v),
\end{equation}
where $ L := \frac{\d {\cal F}}{\d u}( \bar u)$ and
\begin{equation}
\label{g}
\CalG( v)=
g( v, v_x, x):= h(\bar u +v, \bar u_x +v_x) - h(\bar u, \bar u_x) - dh(\bar u, \bar u_x)(v,v_x)
\end{equation}
is a quadratic-order Taylor remainder.

Choosing $\d_t\alpha$ so as to cancel $\Pi_2$ of the righthand
side of \eqref{epspert}, we obtain finally the {\it reduced equations}
\be\label{redv}
\d_t v= \Pi_1( Lv + \CalG(v))
\ee
and
\be\label{alphaeq}
\d_t \alpha =
\frac{\pi_2(Lv + \CalG(v))}
{1+ \pi_2 (\d_x v)}
\ee
for $v\in \phi^\perp$, 
where $\pi_2 v:= \langle \tilde \phi, v\rangle\|\phi|_{L^2}^2$,
of the same regularity as the original equations.

Clearly, \eqref{alphaeq} is well-defined so long as 
$|\d_x v|_{L^\infty}\le C|v|_{H^2}$
remains small, hence we may solve the $v$ equation independently
of $\alpha$, determining $\alpha$-behavior afterward to determine
the full solution 
$$
u(x,t)= \bar u(x-\alpha(t))+v(x-\alpha(t),t).
$$
Moreover, it is easily seen (by the block-triangular structure
of $L$ with respect to this decomposition) that
the linear part $\Pi_1L=\Pi_1 L\Pi_1$ of the $v$-equation possesses
all spectrum of $L$ apart from the zero eigenvalue associated with 
eigenfunction $\phi$.
Thus, we have effectively projected out this zero-eigenfunction, and
with it the group symmetry of translation.

We may therefore construct the center stable manifold for the reduced
equation \eqref{redv}, automatically obtaining translation-invariance
when we extend to the full evolution using \eqref{alphaeq}.
See \cite{TZ1} for further discussion.

\subsection{Preliminary estimates}
For ease of notation, introduce $L_0:=\Pi_1 L$, $\cG_0:=\Pi_1 \cG$.

\begin{lemma}\label{Frechet}
Under the assumed regularity $h\in C^{k+1}$,
both $\CalG$ and $\CalG_0$ are Frechet differentiable of order $(k+1)$
considered respectively as functions from $H^2$ to $H^1$ and
$\phi^\perp\subset H^2$ to $H^1$:
$\CalG$ on the whole space and $\CalG_0$ for $|v|_{H^2}$
sufficiently small.
\end{lemma}

\begin{proof}
Differentiability of $\CalG$ follows by direct calculation;
see \cite{Sat,TZ1}.
Differentiability of $\Pi_1\CalG$ follows similarly, using also
the fact, already discussed, that ${1+ \pi_2 (\d_x v)}$
remains bounded from zero for $|v|_{H^2}$ small, and
the fact (see Lemma \ref{Pi}) that
$\Pi_j$ as bounded linear operators from each $H^s$ to itself
are infinitely differentiable in the Frechet sense.
\end{proof}

\begin{lemma}\label{semigp}
$L_0$ generates an analytic semigroup $e^{ L_0 t}= \Pi_1 e^{Lt}\Pi_1$
on $\phi^\perp \subset H^2$.
Moreover, the unstable (positive real part) spectra of $L$ and $L_0$
agree in both location and multiplicity, with associated total unstable
eigenprojections $\Pi^0_u$ and $\Pi_u$ related by $\Pi_u^0=\Pi_1 \Pi_u \Pi_1$
and total center stable eigenprojections $\Pi^0_{cs}$ and $\Pi_{cs}$ 
related by $\Pi_{cs}^0=\Pi_1 \Pi_{cs} \Pi_1$.
Likewise, except possibly at $\lambda=0$,
the resolvent sets of $L$ and $L_0$ agree, with
$(\lambda-L_0)^{-1}= \Pi_1(\lambda-L)^{-1}\Pi_1$.
\end{lemma}

\begin{proof}
Direct computation using relations $L\Pi_1=L$ and $L\Pi_2=0$
yields the resolvent relation, whence we obtain the remaining
relations by their characterizations in terms of the resolvent
(for example, the characterization of eigenprojection as 
residue of the resolvent operator \cite{Kat}).
As $L$ is a sectorial operator, it follows that $L_0$ is as
well, and both generate analytic semigroups given by the inverse
Laplace transform of the resolvent.
\end{proof}

\begin{corollary}[\cite{TZ1}]\label{linest}
Under assumptions (A0)--(A2),
 \begin{equation} \label{bound-pde} 
\begin{aligned}
\| e^{t L_0}\Pi_{cs}  \|_{H^1\to H^2}  &\leq  C_\o  (1+ t^{-1/2}) e^{\o t}, \\
\| e^{-t L_0}\Pi_u  \|_{H^1\to H^3}  &\leq  C_\o e^{-\b t},\\
\end{aligned}
 \end{equation}
 for some $\b > 0,$ and for all $\o > 0,$ for all $t\ge 0$.
\end{corollary}

\begin{proof}
These follow from the corresponding estimates for $L$, which
are standard semigroup estimates for second-order elliptic operators;
see Appendix \ref{extraproofs}.
\end{proof}

  Let $\rho$ be a smooth truncation function as in \S \ref{lcm} and 
 $ {\cal G}_0^\delta (v) := 
\rho\Big( \frac{ | v |_{H^2}}{\delta}\Big) {\cal G}_0(v).$

 \begin{lemma}[\cite{TZ1}] \label{trunc} The map ${\cal G}^\delta_0: H^2 \times \R \to H^1 \times \R$
  is $C^{k+1}$ and its Lipschitz norm with respect to $v$ 
%can be made arbitrarily small 
is $O(\delta)$ as $\delta\to 0.$
 \end{lemma}

 \begin{proof}
See Appendix A.
\end{proof}

\begin{corollary}\label{integrand}
Under assumptions (A0)--(A2),
 \begin{equation} \label{bound-integrand} 
\begin{aligned}
\| e^{t L_0}\Pi_{cs} {\cal G}^\delta_0
  \|_{H^2\to H^2}  &\leq  
C_\o(1+ t^{-1/2}) e^{\o t}, \\
\| e^{-t L_0}\Pi_u {\cal G}^\delta_0
  \|_{H^2\to H^2}  &\leq  C_\o e^{-\b t}, \\
\end{aligned}
 \end{equation}
for some $\b > 0,$ and for all $\o > 0,$ for all $t\ge 0$,
with Lipshitz bounds
 \begin{equation} \label{bound-integrand-lip} 
\begin{aligned}
\| e^{t L_0}\Pi_{cs} d{\cal G}^\delta_0
  \|_{H^2\to H^2}  &\leq  
C_\o \delta(1+ t^{-1/2}) e^{\o t}, \\
\| e^{-t L_0}\Pi_u d{\cal G}^\delta_0
  \|_{H^2\to H^2}  &\leq  C_\o  \delta e^{-\b t}. \\
\end{aligned}
 \end{equation}
\end{corollary}

%%%%%%%%%%%%%%%%%%%%%%%%%%%%%%%%%%%%%%%%%%%%%%%%%%%%%%%%%%%%%%%%%%%%%
%\subsection{Existence of center stable manifold}\label{pdecsproof}
\subsection{Translation-invariant center stable manifold}\label{pdecsproof}

\begin{proof}[Proof of Proposition \ref{t:cspde}.] 
Using bounds \eqref{bound-pde}, \eqref{bound-integrand}, 
and \eqref{bound-integrand-lip},
and observing that, since finite-dimensional, the unstable flow
$e^{L_0t}\Pi_u$ is well-defined in both forward and backward time,
we find, applying the finite-dimensional argument word-for-word,
that the center stable manifold of the truncated flow
$u_t=L_0u+\cG^\delta_0(u)$ is the graph of the map 
$\Phi$ defined on $\tilde \Sigma_{cs}$ as 
$\Phi(u_{cs}) = (v_{cs}(0)),$ where $v_{cs}$ is the unique solution 
in $\phi^\perp\subset H^2$ of 
 \begin{equation} 
\nonumber
\begin{aligned} v (x,t) = & \, e^{t L_0} u_{cs} 
+ \int_0^t e^{(t - s) L_0} \Pi_{cs} {\cal G}_0^\delta(v(x,s)) ds 
%\\ & 
- \int_t^\infty e^{(t - s) L_0} \Pi_u {\cal G}_0^\delta(v(x,s)) ds
\end{aligned}
\end{equation}
as guaranteed by the Contraction Mapping Theorem.

Indeed, the only deviation from the finite-dimensional case
is the appearance of the new factor
$(1+(t-s)^{-1/2})$ in the integrand of estimates having to do with integral term
$\int_0^te^{L_0(t-s)}\Pi_{cs}\cG^\delta_0 $,
which, since integrable, does not alter the final estimates.
The proof of smoothness relies on these same bounds,
so likewise carries over word-for-word as in the finite-dimensional case.
Finally, invariance in forward time follows by the characterization
of the center stable manifold as the set of solutions of the
truncated equations growing no
faster than $|v(t)|\le Ce^{\beta t}$ in forward time.
(Since the flow of \eqref{rcd} is only a semigroup, 
we cannot conclude invariance in backward time as 
in the finite-dimensional case.)
Since sufficiently small bounded solutions of the original
are bounded solutions also of the truncated equations, they are 
%necessarily 
contained in the center stable manifold, independent
of the choice of truncation function.

This gives a center stable manifold with the stated properties
for the reduced equation \eqref{redv} for $v\in \phi^\perp$.
Solving for the shift $\alpha$ in terms of $v$ using \eqref{alphaeq},
and substituting $v$, $\alpha$ into \eqref{shift} to obtain
$u$, we obtain a translation-invariant center stable manifold
for the original equations \eqref{rcd}.
Observing that small bounded solutions $v$ of the
reduced equations, by \eqref{shift}, 
correspond to solutions $u$ of the \eqref{rcd} 
remaining close to a translate of $\bar u$,
we are done. 
\end{proof}

\subsection{Application to viscous shock waves}\label{viscapp}

\begin{proof} [Proof of Theorem \ref{t:maincs}]
By Proposition \ref{t:cspde}, it is sufficient to verify that
(H0)--(H1) imply (A0)--(A2), for $h(u,u_x):= f(u)_x=df(u)u_x$.
Clearly, (H0) implies (A0) by the form of $h$.
Plugging $u=\bar u(x)$ into \eqref{cons}, we obtain the 
standing-wave ODE $f(\bar u)_x=\bar u_xx$, or, integrating
from $-\infty$ to $x$, the first-order system
$$
\bar u_x= f(\bar u)-f(u_-).
$$
Linearizing about the assumed critical points $u_\pm$ yields
linearized systems $w_t=df(u_\pm)w$, from which we see that
$u_\pm$ are nondegenerate rest points by (H1),
as a consequence of which (A2) follows by standard ODE theory \cite{Co}.

Finally, linearizing PDE \eqref{rcd} about the constant solutions
$u\equiv u_\pm$, we obtain $w_t=L_\pm w:= -df(u_\pm) w_x-w_{xx}$.
By Fourier transform, the limiting operators $L_\pm$ have
spectra $\lambda^\pm_j(k)=-ika^\pm_j(k)-k^2$, where the Fourier
wave-number $k$ runs over all of $\R$; in particular, $L_\pm$
have spectra of nonpositive real part.
By a standard result of Henry \cite{He}, the essential spectrum
of $L$ lies to the left of the rightmost boundary of the spectra
of $L_\pm$, hence we may conclude that
the essential spectrum of $L$ is entirely nonpositive.
As the spectra of $L$ to the right of the essential spectrum
by sectoriality of $L$, consists of finitely
many discrete eigenvalues, this means that the spectra of $L$
with positive real part consists of $p$ unstable eigenvalues,
for some $p$, verifying (A1).
\end{proof}

\section{Conditional stability analysis}\label{s:cond}

Define similarly as in Section \ref{redeq} the perturbation variable
\be\label{pert2}
v(x,t):=u(x+\alpha(t),t)-\bar u(x)
\ee
for $u$ a solution of \eqref{cons}, where $\alpha$ is to be specified later
in a way appropriate for the task at hand.
Subtracting the equations for $u(x+\alpha(t), t)$ and $\bar u(x)$,
we obtain the nonlinear perturbation equation
\be\label{nlpert}
v_t-Lv= N(v)_x,
\ee
where $L:=-\partial_x df(\bar u) +\partial_x^2$ as in \eqref{lineq}
denotes the linearized operator about $\bar u$ and
\be\label{N}
N(v):=-(f(\bar u+ v)-f(\bar u)- df(\bar u)v)
\ee
where, so long as $|v|_{H^1}$ (hence $|v|_{L^\infty}$ and $|u|_{L^\infty}$) 
remains bounded, 
\ba\label{Nbds}
N(v)&=O(|v|^2).\\
%\partial_x N(v)&=O(|v||\partial_x v|),\\
%\partial_x^2 N(v)&= O(|\partial_x|^2+ |v||\partial_x^2v|).\\
\ea

\subsection{Projector bounds}\label{projbds}
Let $\Pi_u$ denote the eigenprojection of $L$ onto its
unstable subspace $\Sigma_u$, and $\Pi_{cs}={\rm Id}- \Pi_u$
the eigenprojection onto its center stable subspace $\Sigma_{cs}$. 

\bl\label{projlem}
Assuming (H0)--(H1),
\be\label{comm}
\Pi_j \partial_x= \partial_x \tilde \Pi_j
\ee
for $j=u,\, cs$ and, for all $1\le p\le \infty$, $0\le r\le 4$, 
\ba
|\Pi_{cs}|_{W^{r,p}\to W^{r,p}}, |\tilde \Pi_{cs}|_{W^{r,p}\to W^{r,p}}&\le C,\\
|\tilde \Pi_{cs}|_{W{r,p}\to W^{r,p} }, 
\; |\tilde \Pi_{cs}|_{W{r,p}\to W^{r,p}}&\le C.
\ea
\el

\begin{proof}
Recalling (see the proof of Theorem \ref{t:maincs})
that $L$ has at most finitely many unstable eigenvalues,
we find that $\Pi_u$ may be expressed as
$$
\Pi_u f= \sum_{j=1}^p \phi_j(x) \langle \tilde \phi_j, f\rangle,
$$
where $\phi_j$, $j=1, \dots p$ are generalized right eigenfunctions of
$L$ associated with unstable eigenvalues $\lambda_j$, 
satisfying the generalized eigenvalue equation $(L-\lambda_j)^{r_j}\phi_j=0$,
$r_j\ge 1$, and $\tilde \phi_j$
are generalized left eigenfunctions.
Noting that $L$ is divergence form, and that $\lambda_j\ne 0$,
we may integrate $(L-\lambda_j)^{r_j}\phi_j=0$ over $\R$ to
obtain $\lambda_j^{r_j}\int \phi_j dx=0$ and thus $\int\phi_j dx=0$.
Noting that $\phi_j$, $\tilde \phi_j$ and derivatives decay exponentially
by standard theory \cite{He,ZH,MaZ1}, we find that
$$
\phi_j= \partial_x \Phi_j
$$
with $\Phi_j$ and derivatives exponentially decaying, hence
$$
\tilde \Pi_u f=\sum_j \Phi_j \langle \partial_x \tilde \phi, f\rangle.
$$
Estimating 
$$
|\partial_x^j\Pi_u f|_{L^p}=|\sum_j \partial_x^j\phi_j \langle \tilde \phi_j f
\rangle|_{L^p}\le
\sum_j |\partial_x^j\phi_j|_{L^p} |\tilde \phi_j|_{L^q} |f|_{L^p}
\le C|f|_{L^p}
$$
for $1/p+1/q=1$
and similarly for $\partial_x^r \tilde \Pi_u f$, we obtain the claimed
bounds on $\Pi_u$ and $\tilde \Pi_u$, from which the bounds on
$\Pi_{cs}={\rm Id}-\Pi_u$ and
$\tilde \Pi_{cs}={\rm Id}-\tilde \Pi_u$ follow immediately.
\end{proof}

\subsection{Linear estimates}

Let $G_{cs}(x,t;y):=\Pi_{cs}e^{Lt}\delta_y(x)$ denote
the Green kernel of the linearized solution operator on
the center stable subspace $\Sigma_{cs}$.
Then, we have the following detailed pointwise bounds
established in \cite{TZ2,MaZ1}.

\begin{proposition}[\cite{TZ2,MaZ1}]\label{greenbounds}
Assuming (H0)--(H2), (D1)--D(3), 
the center stable Green function may be decomposed as 
$G_{cs}=E+\tilde G$, where 
\begin{equation}\label{E}
E(x,t;y)= \d_x \bar u(x) e_j(y,t),
\end{equation}
\begin{equation}\label{e}
  e(y,t)=\sum_{a_k^{-}>0}
  \left(\textrm{errfn }\left(\frac{y+a_k^{-}t}{\sqrt{4t}}\right)
  -\textrm{errfn }\left(\frac{y-a_k^{-}t}{\sqrt{4t}}\right)\right)
  l_{k}^{-}(y)
\end{equation}
for $y\le 0$ and symmetrically for $y\ge 0$, 
$l_k^-\in \R^n$ constant, and
%NOTE: e_j continuuous at zero, though we don't need it....
\begin{equation}\label{Gbounds}
\begin{aligned}
|\tilde G(x,t;y)|&\le  Ce^{-\eta(|x-y|+t)} +
\sum_{k=1}^n t^{-1/2}e^{-(x-y-a_k^{-} t)^2/Mt} e^{-\eta x^+} \\
&+
\sum_{a_k^{-} > 0, \, a_j^{-} < 0} 
\chi_{\{ |a_k^{-} t|\ge |y| \}}
t^{-1/2} e^{-(x-a_j^{-}(t-|y/a_k^{-}|))^2/Mt}
e^{-\eta x^+} \\
&+
\sum_{a_k^{-} > 0, \, a_j^{+}> 0} 
\chi_{\{ |a_k^{-} t|\ge |y| \}}
t^{-1/2} e^{-(x-a_j^{+} (t-|y/a_k^{-}|))^2/Mt}
e^{-\eta x^-}, \\
\end{aligned}
\end{equation}
\begin{equation}\label{Gybounds}
\begin{aligned}
|\partial_y \tilde G(x,t;y)|&\le  Ce^{-\eta(|x-y|+t)}
+ Ct^{-1/2}\Big( \sum_{k=1}^n 
t^{-1/2}e^{-(x-y-a_k^{-} t)^2/Mt} e^{-\eta x^+} \\
&+
\sum_{a_k^{-} > 0, \, a_j^{-} < 0} 
\chi_{\{ |a_k^{-} t|\ge |y| \}}
t^{-1/2} e^{-(x-a_j^{-}(t-|y/a_k^{-}|))^2/Mt}
e^{-\eta x^+} \\
&+
\sum_{a_k^{-} > 0, \, a_j^{+}> 0} 
\chi_{\{ |a_k^{-} t|\ge |y| \}}
t^{-1/2} e^{-(x-a_j^{+} (t-|y/a_k^{-}|))^2/Mt}
e^{-\eta x^-}\Big) \\
\end{aligned}
\end{equation}
for $y\le 0$ and symmetrically for $y\ge 0$,
for some $\eta$, $C$, $M>0$, where 
$a_j^\pm$ are the eigenvalues of $df(u_\pm)$, 
$x^\pm$ denotes the positive/negative
part of $x$, and  indicator function $\chi_{\{ |a_k^{-}t|\ge |y| \}}$ is 
$1$ for $|a_k^{-}t|\ge |y|$ and $0$ otherwise.
\end{proposition}

\begin{proof}
As observed in \cite{TZ2},
it is equivalent to establish decomposition 
\be\label{fulldecomp}
G=G_u + E+\tilde G
\ee
for the full Green function $G(x,t;y):=e^{Lt}\delta_y(x)$,
where 
$$
G_u(x,t;y):=\Pi_u e^{Lt}\delta_y(x)
=
e^{\gamma t}\sum_{j=1}^p\phi_j(x)\tilde \phi_j(y)^t
$$
for some constant matrix $M\in \C^{p\times p}$
denotes the Green kernel of the linearized solution operator
on $\Sigma_u$, $\phi_j$ and $\tilde\phi_j$ right and left
generalized eigenfunctions associated with unstable eigenvalues
$\lambda_j$, $j=1,\dots,p$.

The problem of describing the full Green function
has been treated in \cite{ZH, MaZ3}, 
starting with the Inverse Laplace Transform representation
\be\label{ILT2}
G(x,t;y)=e^{Lt}\delta_y(x)= \oint_\Gamma e^{\lambda t}(\lambda-L(\e))^{-1} 
\delta_y(x)d\lambda \, ,
\ee
where 
$$
\Gamma:= \partial \{ \lambda : \Re \lambda\le \eta_1 - \eta_2 |\Im \lambda|\}
$$
is an appropriate sectorial contour, $\eta_1$, $\eta_2>0$;
estimating the resolvent kernel 
$G^\eps_\lambda(x,y):=(\lambda-L(\e))^{-1}\delta_y(x)$
using Taylor expansion in $\lambda$,
asymptotic ODE techniques in $x$, $y$, and judicious decomposition
into various scattering, excited, and residual modes;
then, finally, estimating the contribution of various modes to \eqref{ILT2}
by Riemann saddlepoint (Stationary Phase) method, moving contour
$\Gamma$ to a optimal, ``minimax'' positions for each
mode, depending on the values of $(x,y,t)$.

In the present case, we may first move $\Gamma$ to a contour
$\Gamma'$ enclosing (to the left) all spectra of $L$
except for the $p$ unstable eigenvalues $\lambda_j$, $j=1, \dots, p$,
to obtain
$$
G(x,t;y)= \oint_{\Gamma'} e^{\lambda t}(\lambda-L)^{-1} d\lambda
+ \sum_{j=\pm} 
\Res_{\lambda_j(\eps)} \big( e^{\lambda t}(\lambda-L)^{-1}
\delta_y(x) \big),
$$
where
$\Res_{\lambda_j(\eps)} \big( e^{\lambda t}(\lambda-L)^{-1}
\delta_y(x) \big) = G_u(x,t;y)$, then estimate the remaining term
$ \oint_{\Gamma'} e^{\lambda t}(\lambda-L)^{-1} d\lambda$
on minimax contours as just described.
See the proof of Proposition 7.1, \cite{MaZ3}, for a detailed 
discussion of minimax estimates $E+G$ and of Proposition 7.7, 
\cite{MaZ3} for a complementary discussion of residues
incurred at eigenvalues in $\{\Re \lambda\ge 0\}\setminus\{0\}$.
See also \cite{TZ1}.
\end{proof}

\bc [\cite{MaZ1}] \label{lpbds}
Assuming (H0)--(H2), (D1)--(D3),
\be \label{tGbounds}
|\int_{-\infty}^{+\infty} \tG(\cdot,t;y)f(y)dy|_{L^p}
\le C (1+t)^{-\frac{1}{2}(\frac{1}{q}-\frac{1}{p})} |f|_{L^q},
\ee
\be\label{tGybounds}
|\int_{-\infty}^{+\infty} \tG_y(\cdot,t;y)f(y)dy|_{L^p}
\le C (1+t)^{-\frac{1}{2}(\frac{1}{q}-\frac{1}{p})-\frac{1}{2}} |f|_{L^q},
\ee
for all $t\ge 0$, some $C>0$, for any
$1\le q\le p$ (equivalently, $1\le r\le p$)
and $f\in L^q$, where $1/r+1/q=1+1/p$.
\ec

\begin{proof}
Standard convolution inequalities together
with bounds \eqref{Gbounds}--\eqref{Gybounds}; see \cite{MaZ1,MaZ2,MaZ3,Z2} 
for further details.
\end{proof}

\bc[\cite{Z4}]\label{ebds}
The kernel ${e}$ satisfies
$$
|{e}_y (\cdot, t)|_{L^p},  |{e}_t(\cdot, t)|_{L^p} 
\le C t^{-\frac{1}{2}(1-1/p)},
\label{36}
$$
$$
|{e}_{ty}(\cdot, t)|_{L^p} 
\le C t^{-\frac{1}{2}(1-1/p)-1/2},
\label{37}
$$
for all $t>0$.  
%%%%%
%BEGIN CHANGE, March 24
Moreover, for $y\le 0$ we have the pointwise bounds
$$
|{e}_y (y,t)|, |{e}_t (y,t)| 
\le 
Ct^{-\frac{1}{2}} \sum_{a_k^->0}
\Big(e^{-\frac{(y+a_k^-t)^2}{Mt}}+
e^{-\frac{(y-a_k^-t)^2}{Mt}}\Big),
\label{38}
$$
$$
|{e}_{ty} (y,t)| \le 
Ct^{-1} \sum_{a_k^->0}
\Big(e^{-\frac{(y+a_k^-t)^2}{Mt}}+
e^{-\frac{(y-a_k^-t)^2}{Mt}}\Big),
\label{39}
$$
for $M>0$ sufficiently large, 
and symmetrically for $y\ge 0$.
\ec

\begin{proof}
Direct computation using with definition \eqref{e}; 
see Appendix \ref{extraproofs}.
\end{proof}

%NOTE: to get higher $x$-derivs, differentiate equations and keep
%principal part plus lower-order commutator terms, proceed by iteration.
%details ommitted... (no problem here in fact...)

\subsection{Reduced equations II}
%(instantaneous projection)
Recalling that $\d_x\bar u$ is a stationary
solution of the linearized equations $u_t=Lu$,
so that $L\d_x\bar =0$, or
$$
\int^\infty_{-\infty}G(x,t;y)\bar u_x(y)dy=e^{Lt}\bar u_x(x)
=\d_x\bar u(x),
$$
we have, applying Duhamel's principle to \eqref{nlpert},
$$
\begin{array}{l}
  \displaystyle{
  v(x,t)=\int^\infty_{-\infty}G(x,t;y)v_0(y)\,dy } \\
  \displaystyle{\qquad
  -\int^t_0 \int^\infty_{-\infty} G_y(x,t-s;y)
  (N(v)+\dot \alpha v ) (y,s)\,dy\,ds + \alpha (t)\d_x \bar u(x).}
\end{array}
$$
Defining 
\begin{equation}
 \begin{array}{l}
  \displaystyle{
  \alpha (t)=-\int^\infty_{-\infty}e(y,t) v_0(y)\,dy }\\
  \displaystyle{\qquad
  +\int^t_0\int^{+\infty}_{-\infty} e_{y}(y,t-s)(N(v)+
  \dot \alpha\, v)(y,s) dy ds, }
  \end{array}
 \label{alpha}
\end{equation}
following \cite{ZH,Z4,MaZ2,MaZ3}, 
where $e$ is defined as in \eqref{e}, 
and recalling the decomposition $G=E+ G_u+ \tilde G$ of \eqref{fulldecomp},
we obtain the {\it reduced equations}
\begin{equation}
\begin{array}{l}
 \displaystyle{
  v(x,t)=\int^\infty_{-\infty} (G_u+\tilde G)(x,t;y)v_0(y)\,dy }\\
 \displaystyle{\qquad
  -\int^t_0\int^\infty_{-\infty}(G_u+\tilde G)_y(x,t-s;y)(N(v)+
  \dot \alpha v)(y,s) dy \, ds, }
\end{array}
\label{v}
\end{equation}
and, differentiating (\ref{alpha}) with respect to $t$,
and observing that 
$e_y (y,s)\rightharpoondown 0$ as $s \to 0$, as the difference of 
approaching heat kernels,
\begin{equation}
 \begin{array}{l}
 \displaystyle{
  \dot \alpha (t)=-\int^\infty_{-\infty}e_t(y,t) v_0(y)\,dy }\\
 \displaystyle{\qquad
  +\int^t_0\int^{+\infty}_{-\infty} e_{yt}(y,t-s)(N(v)+
  \dot \alpha v)(y,s)\,dy\,ds. }
 \end{array}
\label{alphadot}
\end{equation}
\medskip

We emphasize that this (nonlocal in time) choice of $\alpha$ 
and the resulting reduced equations are completely different from those
of Section \ref{redeq}, according to their respective purposes.
A third possible choice has been introduced in \cite{GMWZ1,GMWZ2} 
for the study of the inviscid limit problem.
As discussed further in \cite{Go,Z4,MaZ2,MaZ3,GMWZ1,BeSZ},
$\alpha$ may be considered in the present context as defining a notion of
approximate shock location.

\subsection{Nonlinear damping estimate}

\begin{proposition}[\cite{MaZ3}]\label{damping}
Assuming (H0)-(H3), let $v_0\in H^{4}$, 
%$k\ge 2$ as in (H0), 
and suppose that for $0\le t\le T$, the $H^{4}$ norm of $v$
remains bounded by a sufficiently small constant, for $v$ as in
\eqref{pert2} and $u$ a solution of \eqref{cons}.
Then, for some constants $\theta_{1,2}>0$, for all $0\leq t\leq T$,
\begin{equation}\label{Ebounds}
\|v(t)\|_{H^4}^2 \leq C e^{-\theta_1 t} \|v(0)\|^2_{H^4} 
+ C \int_0^t e^{-\theta_2(t-s)} (|v|_{L^2}^2 + |\dot \alpha|^2) (s)\,ds.
\end{equation}
\end{proposition}

\begin{proof}
Subtracting the equations for $u(x+\alpha(t), t)$ and $\bar u(x)$,
we may write the perturbation equation for $v$ alternatively as
\begin{equation}\label{vperturteq}
v_t + \left( \int_0^1 df \big(\bar u(x) + \tau v(x,t) \big)\,d\tau\, 
v\right)_x - v_{xx} =
\dot{\alpha}(t) \d_x \bar u(x).
\end{equation}
Observing that $\partial_x^j(\partial_x \bar{u})(x)=O(e^{-\eta|x|})$ 
is bounded in $L^1$ norm for $j\le 4$, we take the $L^2$ inner product 
in $x$ of $\sum_{j=0}^4\partial_x^{2j}v$ against \eqref{vperturteq}, 
integrate by parts and rearrange the resulting terms to arrive 
at the inequality
\[
\partial_t \|v\|_{H^4}^2(t) \leq -\theta \|\partial_x^{5} v\|_{L^2}^2 +
C \left(\|v\|_{H^4}^2 + |\dot{\alpha}(t)|^2\right),
\]
$\theta>0$, for $C>0$ sufficiently large, so long as $\|v\|_{H^4}$ 
remains bounded. Using the Sobolev interpolation
\[
\|v\|_{H^4}^2 \leq \tilde{C}^{-1} \|\partial_x^{5} v\|_{L^2}^2 + \tilde{C} \| v\|_{L^2}^2
\]
for $\tilde{C}>0$ sufficiently large, we obtain 
\[
\partial_t \|v\|_{H^4}^2(t) \leq -\tilde{\theta} \|v\|_{H^4}^2 + C \left(\|v\|_{L^2}^2 + |\dot{\alpha}(t)|^2\right),
\]
from which \eqref{Ebounds} follows by Gronwall's inequality.
\end{proof}

\subsection{Proof of nonlinear stability}

Decompose now the nonlinear perturbation $v$ as
\be\label{vdecomp}
v(x,t)=w(x,t)+z(x,t),
\ee
where
\be\label{wzdef}
w:=\Pi_{cs}v, \quad z:=\Pi_u v.
\ee
Applying $\Pi_{cs}$ to \eqref{v} and recalling commutator
relation \eqref{comm}, we obtain an equation
\ba \label{w}
  w(x,t)&=\int^\infty_{-\infty} \tilde G (x,t;y)w_0(y)\,dy \\
  &\quad -\int^t_0\int^\infty_{-\infty} \tilde G_y (x,t-s;y)
\tilde \Pi_{cs}(N(v)+
  \dot \alpha v)(y,s) dy \, ds
\ea
for the flow along the center stable manifold, parametrized by
$w\in \Sigma_{cs}$.

\bl\label{quadlem} Assuming (H0)--(H1), for $v$ lying initially
on the center stable manifold $\cM_{cs}$,
\be\label{vwbd}
|z|_{W^{r,p}}\le C|w|_{H^2}^2
\ee 
for some $C>0$, for all $1\le p\le \infty$ and $0\le r\le 4$,
so long as $|w|_{H^2}$ remains sufficiently small.
\el

\begin{proof}
By tangency of the center stable manifold to $\Sigma_{cs}$, we
have immediately $|z|_{H^2}\le C|w|_{H^2}^2$, whence
\eqref{vwbd} follows by equivalence of norms for finite-dimensional
vector spaces, applied to the $p$-dimensional subspace $\Sigma_u$.
(Alternatively, we may see this by direct computation using
the explicit description of $\Pi_u v$ afforded by Lemma \ref{projlem}.)
\end{proof}

\begin{proof}[Proof of Theorem \ref{t:mainstab}]

%instab. part:
Recalling by Theorem \ref{t:maincs}
that solutions remaining for all time in a sufficiently
small radius neighborhood $\cN$ of the set of translates of $\bar u$
lie in the center stable manifold $\cM_{cs}$, we obtain trivially
that solutions not originating in $\cM_{cs}$ must exit $\cN$ in finite time,
verifying the final assertion of orbital instability with respect
to perturbations not in $\cM_{cs}$.

Consider now a solution $v \in \cM_{cs}$, or, equivalently,
a solution $w\in \Sigma_{cs}$ of \eqref{w} with 
$z=\Phi_{cs}(w)\in \Sigma_u$.
Define
\begin{equation}
\label{zeta2}
 \zeta(t):= \sup_{0\le s \le t}
 \Big( |w|_{H^2}(1+s)^{\frac{1}{4}} + |w|_{L^\infty}+
|\dot \alpha (s)|(1+s)^{\frac{1}{2}} \Big).
\end{equation}
We shall establish:

{\it Claim.} For all $t\ge 0$ for which a solution exists with
$\zeta$ uniformly bounded by some fixed, sufficiently small constant,
there holds
\begin{equation}
\label{claim}
\zeta(t) \leq C_2(E_0 + \zeta(t)^2)
\quad \hbox{\rm for} \quad
E_0:=|v_0|_{L^1\cap H^2}.
\end{equation}
\medskip

{}From this result, provided $E_0 < 1/4C_2^2$, 
we have that $\zeta(t)\le 2C_2E_0$ implies
$\zeta(t)< 2C_2E_0$, and so we may conclude 
by continuous induction that
 \begin{equation}
 \label{bd}
  \zeta(t) < 2C_2E_0
 \end{equation}
for all $t\geq 0$, whence we obtain the stated bounds by definition
\eqref{zeta2}.
(By Lemma \ref{damping} and standard short-time $H^s$ existence theory, 
$v\in H^4$ exists and $\zeta$ remains
continuous so long as $\zeta$ remains bounded by some uniform constant,
hence \eqref{bd} is an open condition.)
Thus, it remains only to establish the claim above.
\medskip

{\it Proof of Claim.}
We must show that $u(\theta+\psi_1+\psi_2)^{-1}$ and
$|\dot \alpha(s)|(1+s)$
are each bounded by
$C(E_0 + \zeta(t)^2)$,
for some $C>0$, all $0\le s\le t$, so long as $\zeta$ remains
sufficiently small.

By Lemma \ref{quadlem},
$ |w_0|_{L^1\cap H^2}\le |v_0|_{L^1\cap H^2}+ |z_0|_{L^1\cap H^2}
\le
|v_0|_{L^1\cap H^2}+ C|w_0|_{H^2}^2, $
whence
$$
|w_0|_{L^1\cap H^2}\le CE_0.
$$
Likewise, by Lemma \ref{quadlem},
\eqref{zeta2}, \eqref{Nbds}, and Lemma \ref{projlem}, for $0\le s\le t$,
\ba\label{Nlast}
|\tilde \Pi_{cs}(N(v)+ \dot \alpha v)(y,s)|_{L^q}&\le C\zeta(t)^2 
(1+s)^{-\frac{5}{4}}.\\
%|\tilde \Pi_{cs}(N(v)+ \dot \alpha v)(y,s)|_{H^2}&\le C\zeta(t)^2 
%(1+s)^{-\frac{5}{4}}.
\ea

Combining the latter bounds with representations \eqref{w}--\eqref{alphadot},
taking $q=2$ in \eqref{Nlast}, and applying Corollary \ref{lpbds}, we obtain
 \ba\label{claimw}
  |w(x,t)|_{L^p} &\le
  \Big|\int^\infty_{-\infty} \tilde G(x,t;y) w_0(y)\,dy\Big|_{L^p}
   \\
 &\qquad +
\Big|\int^t_0
  \int^\infty_{-\infty} \tilde G_y(x,t-s;y)  \tilde \Pi_{cs}(N(v)+
  \dot \alpha v)(y,s)  dy \, ds\Big|_{L^p} \\
  & \le
  E_0 (1+t)^{-\frac{1}{2}(1-\frac{1}{p})}
   \\
 &\quad +
C\zeta(t)^2 \int^t_0
  \int^\infty_{-\infty}(t-s)^{-\frac{3}{4}+\frac{1}{2p}}
(1+s)^{-\frac{3}{4}} dy \, ds \\
&\le
C(E_0+\zeta(t)^2)(1+t)^{-\frac{1}{2}(1-\frac{1}{p})}
\ea
and, similarly, using H\"older's inequality and
applying Corollary \ref{ebds},
\ba\label{claimalpha}
 |\dot \alpha(t)| &\le \int^\infty_{-\infty}|e_t(y,t)|
  |v_0(y)|\,dy \\
  &\qquad +\int^t_0\int^{+\infty}_{-\infty} |e_{yt}(y,t-s)||\tilde \Pi_{cs}
(N(v)+ \dot \alpha v)(y,s)|\,dy\,ds\\
&\le |e_t|_{L^\infty} |v_0|_{L^1}
+ C\zeta(t)^2 \int^t_0
|e_{yt}|_{L^2}(t-s) |\tilde \Pi_{cs}(N(v)+ \dot \alpha v)|_{L^2}(s) ds\\
&\le E_0 (1+t)^{-\frac{1}{2}}
+ C\zeta(t)^2 \int^t_0
(t-s)^{-\frac{3}{4}}(1+s)^{-\frac{3}{4}} ds\\
&\le C(E_0+\zeta(t)^2)(1+t)^{-\frac{1}{2}}.\\
\ea
Applying Lemma \ref{damping} and using \eqref{claimw} and
\eqref{claimalpha},
we obtain, finally,
\be\label{claimwH2}
|w|_{H^2}(t)\le C(E_0+\zeta(t)^2)(1+t)^{-\frac{1}{4}}.
\ee
Combining \eqref{claimw}, \eqref{claimalpha}, and \eqref{claimwH2},
we obtain \eqref{claim} as claimed, completing the proof of
the Claim and the theorem.
\end{proof}

\br\label{ucpt}
\textup{
We point out that the finite-dimensional part $z$ of $v$
is in fact controlled pointwise by its $L^2$ norm and thus
by $|w|_{H^2}$, satisfying
$ |z(x,t)|\le Ce^{-\theta |x|}|z(\cdot, t)|_{L^2(x)}
\le   Ce^{-\theta |x|}|w(\cdot, t)|_{H^2(x)}^2 , $
making possible a pointwise version of the argument above.
%NOTE: (That is, just as for $L^p$ analysis, $z$ is negligible, can essentially
%be ignored... 
This is a key point in treating the nonclassical 
over- or undercompressive cases, 
which appear to require pointwise bounds \cite{HZ,RZ}.
}
\er

\medbreak
{\bf Acknowledgement.}
Thanks to Milena Stanislavova and Charles Li for two interesting discussions
that inspired this work, and to Milena Stanislavova
for pointing out the reference in \cite{GJLS}.

%%%%%%%%%%%%%%%%%%%

%APPENDIX

%%%%%%%%%%%%%%
%
\appendix

\section{Proofs of miscellaneous lemmas}\label{extraproofs}

We include for completeness the proofs of earlier cited lemmas
%TODO: delete second part?
that were not proved in the main body of the text.

\begin{proof}[Proof of Lemma \ref{poseta}]
By Taylor's Theorem,
$$
g(f_2(t))- g(f_1(t))= dg(f_1(t))(f_2(t)-f_1(t)) + o(|f_2(t)-f_1(t)|),
$$
as $|f_2(t)-f_1(t)|\to 0$.
By the assumed uniform boundedness of $|dg|$, we readily obtain
that $f\to dg f$ is a bounded linear operator from
$\cB_\eta\to \cB_\eta$.
On the other hand $|f_2-f_1|(t)\le e^{-\eta (t-t_0)}\|f_2-f_1\|_\eta$
for $t\ge t_0$ implies $|f_2-f_1|\to 0$ uniformly on $t\ge t_0$
as $\|f_2-f_1\|_\eta\to 0$, and so we have
\ba
\|g(f_2(t))- g(f_1(t))- &dg(f_1(t))(f_2(t)-f_1(t))\|_\eta =\\
&\sup_{t\ge t_0} e^{\eta (t-t_0)}o(|f_2(t)-f_1(t)|) =o(\|f_2-f_1\|_\eta,
\ea
yielding the result for $k=1$.  The general result 
then follows by induction on $k$,
applying the result for $k=1$ to successively higher derivatives of $g$.
\end{proof}

\begin{proof}[Proof of Lemma \ref{fixeddep}]
(i) Triangulating, we have
\ba
|y(x_2)-y(x_1)|&=|T(x_2,y(x_2))- T(x_1,y(x_1)|\\
&\le 
|T(x_2,y(x_2))- T(x_2,y(x_1)|+ |T(x_2,y(x_1))- T(x_1,y(x_1)|\\
&\le 
\theta |y(x_2)- y(x_1)|+ L|x_2- x_1|,
\ea
where $0<\theta<1$ and $0<L$ are contraction and Lipschitz coefficients,
yielding after rearrangement
$ |y(x_2)-y(x_1)|\le \frac{ L}{1-\theta} |x_2- x_1|$.

(ii) Applying Taylor's Theorem, and using the result of (i), we have
\ba
y(x_2)-y(x_1)&=T_x(x_2-x_1) + T_y(y_2-y_1) +o(|x_2-x_1|+|y(x_2)-y(x_1)|)\\
&=
T_x(x_2-x_1) + T_y(y_2-y_1) +o(|x_2-x_1|),
\ea
where all derivatives are evaluated at $(x_1,y(x_1))$.
Noting that the operator norm $|T_y|$ is bounded by contraction 
coefficient $0<\theta<1$, 
we have by Neumann series expansion that
$(\Id -T_y)$ is invertible with uniformly bounded inverse
$|(\Id -T_y)^{-1}|\le (1-\theta)^{-1}$.
Thus, rearranging, we have
$$
y(x_2)-y(x_1)=(\Id-T_y)^{-1}T_x(x_2-x_1)+ o(|x_2-x_1|),
$$
yielding the result for $k=1$ by definition of (Frechet) derivative.
The results for $k\ge 1$ then follow by induction upon differentiation of
\eqref{derform}.
\end{proof}

\begin{proof}[Proof of Lemma \ref{truncbds}]
The Lipshitz constant is bounded by $\max |\partial_w N^\eps|$, where
$$
\begin{aligned}
|\partial_w N^\eps|&= |(\rho^\eps)' N+ \rho^\eps \partial_w N|\\
&= |(|w|/\eps)\rho'(|w|/\eps) (N(t,w)/|w|)+ \rho(|w|/\eps) \partial_w N(t,w)|\\
&\le 2\eps \Big(
\max |r|\rho'(|r|) \max_{|w|\le 2\eps_0}\frac{|N(t,w)|}{|w|^2}
+ \max_{|w|\le 2\eps_0}\frac{|\partial_w N(t,w)|}{|w|}\Big)\\
&\le 2\eps \Big(
\max |r|\rho'(|r|) \max_{|w|\le 2\eps_0}
+ 1 \Big) \max_{|w|\le 2\eps_0}|\partial_w ^2N(t,w)|,
\end{aligned}
$$
the final inequality following by $N(t,0)\equiv 0$, 
$\partial_w N(t,0)\equiv 0$ and
the Integral Mean Value Theorem, or first-order Taylor
remainder formula.
\end{proof}

\begin{proof}[Proof of Corollary \ref{linest}]
By Lemma \ref{semigp}, it is sufficient to prove the corresponding
bounds for the purely differential operator $L$.
By sectoriality of $L$,
we have the inverse Laplace transform representations
\begin{equation}\label{ILT}
\begin{aligned}
 e^{tL}\Pi_u &:= \int_{{\Gamma_u }} e^{\lambda t} 
(\lambda - L)^{-1} \, d\lambda,\\
 e^{tL}\Pi_{cs} &:= \int_{{\Gamma_{cs}}} e^{\lambda t} 
(\lambda - L)^{-1} \, d\lambda,\\
\end{aligned}
\end{equation}
where $\Gamma_{cs}$ denotes a sectorial contour bounding the center and stable
spectrum to the right \cite{Pa}, which by (A1) may be taken
so that $\Re \Gamma_s\le \o$, and $\Gamma_u$ denotes a closed curve enclosing
the unstable spectrum of $L$, with $\Re \Gamma_s\ge \beta>0$.

Applying the resolvent formula 
$ L(\lambda-L)^{-1}= \lambda(\lambda-L)^{-1}- Id $,
we obtain in the standard way
$$
 e^{tL}\Pi_j := \int_{{\Gamma_j}} \lambda e^{\lambda t} 
(\lambda - L)^{-1} \, d\lambda,
$$
from which we obtain immediately the second stated bound, and,
by a scaling argument \cite{Pa}, the bound 
\begin{equation}
\| e^{t L} \Pi_{cs} \|_{H^1\to H^3} 
\leq 
\|L e^{t L} \Pi_s \|_{H^1\to H^1} 
\leq C(1+ t^{-1}) e^{\o t}. 
\end{equation}
Recalling the standard bound
$\| e^{t L} \Pi_{cs} \|_{H^1\to H^1} \leq C e^{\o t}$,
and interpolating between $|\cdot|_{H^1}$ and $|\cdot|_{H^3}$,
we obtain the first stated bound. 
\end{proof}

 \begin{proof}[Proof of Lemma \ref{trunc}]
%NOTE: fixed, was incomplete in TZ1...
 The norm in $H^2$ is a quadratic form, hence the map
 $$ v \in H^2 \mapsto \rho\Big( \frac{ | v |_{H^2}}{\delta}\Big) 
\in \R_+,$$
 is smooth, and ${\cal G}_0^\delta$ is as regular as ${\cal G}_0.$ Now
 \begin{eqnarray} | {\cal G}_0^\delta(v_1) - 
{\cal G}_0^\delta(v_2)|_{H^1} 
& \leq & | \rho\Big( \frac{ | v_1 |_{H^2}}{\delta}\Big) -   
 \rho\Big( \frac{ | v_2 |_{H^2}}{\delta}\Big) |_{L^\infty} 
| {\cal G}_0(v_1)|_{H^1} \nonumber \\ &  & + \, | 
\rho\Big( \frac{ | v_2 |_{H^2}}{\delta}\Big) |_{L^\infty}  
| {\cal G}_0(v_1) - {\cal G}_0(v_2)|_{H^1} \nonumber \\
 & \leq & 3 | v_1 - v_2|_{H^2} 
\Big(
\sup_{|v|_{H^2} < \delta} \frac{| {\cal G}_0(v)|_{H^1}}{\delta}
+   
\sup_{|v|_{H^2} < \delta} | d{\cal G}_0(v)|_{H^1}
\Big) ,\nonumber
 \end{eqnarray}
 and $\sup_{|v|_{H^2} < \delta}  | {\cal G}_0(v)|_{H^1} = O( \delta^2),$
 $\sup_{|v|_{H^2} < \delta}  | d{\cal G}_0(v)|_{H^1} = O( \delta).$
 \end{proof}

\begin{proof}[Proof of Corollary \ref{ebds}]
For definiteness, take $y \le 0$. 
Then, \eqref{e} gives
$$
e_y(y,t):= \sum_{a_k^{-}>0} [c_{k,-}^0]l_k^{-t}
\left(K(y+a_k^{-}t, t)
-K (y-a_k^{-}t,t)\right)
\label{41}
$$
$$
e_t(y,t):= \sum_{a_k^{-}>0} [c_{k,-}^0]l_k^{-t}
  \left( (K+K_y)(y+a_k^-t,t)-(K+K_y)(y-a_-k^t,t)\right),
\label{41a}
$$
$$
e_{ty}(y,t):= \sum_{a_k^{-}>0} [c_{k,-}^0]l_k^{-t}
  \left( (K_y+K_{yy})(y+a_k^-t,t)-(K_y+K_{yy})(y-a_k^-t,t)\right),
\label{42}
$$
where
$$
 K(y,t):= \frac{e^{-y^2/4t}}{\sqrt{4\pi t}}
\label{43}
$$
denotes the standard heat kernel.
The pointwise bounds \eqref{38}--\eqref{39} follow immediately for $t\ge 1$ 
by properties of the heat kernel, in turn yielding \eqref{36}--\eqref{37} 
in this case.  
The bounds for small time $t\le 1$ follow from estimates
$$
\aligned
 |K_y (y+at,t, \beta)-K_y (y-at,t, \beta)|
 &=\left|\int^{y-at}_{y+at}K_{yy}(z,t, \beta)\,dz\right|  \\
  & \le Ct^{-3/2}\int^{y-at}_{y+at} e^{\frac{-z^2}{Mt}}\,dz
  \le Ct^{-1/2}e^{-\frac{y^2}{Mt}},
\endaligned
\label{44}
$$
and, similarly,
 $$
\aligned
   |K_{yy}(y+at,t, \beta)-K_{yy}(y-a,t, \beta)|
   &=\left|\int^{at}_{-at}K_{yyy}(z,t, \beta)\,dz\right|
    \\
  & \le Ct^{-2}\int^{y-at}_{y+at} e^{\frac{-z^2}{Mt}}\,dz
   \le Ct^{-1}e^{-\frac{y^2}{Mt}}.
\endaligned
\label{45}
 $$
The bounds for $|{e}_y|$ are again immediate.
\end{proof}

%%%%%%%%%%%%%%

\end{document}